\documentclass[12pt]{article}
\usepackage{amssymb,amsmath,amsthm}
\usepackage{bm}

\def\suml{\mathop{\sum}\limits}                                 
\def\maxl{\mathop{\max}\limits}                                 
\def\e{\varepsilon}                                             
\def\G{\Gamma}                                                  
\def\cd{\!\cdot}                                                
\DeclareMathSymbol{\ell}{\mathord}{letters}{96}                 
\def\id{\operatorname{id}}                                      
\def\od{\operatorname{od}}                                      
\def\FF{\mathop{\cal F}\nolimits}                               
\def\RR{\mathop{\cal R}\nolimits}                               
\def\PP{\mathop{\cal P}\nolimits}                               
\def\rr{\mathop{R}\nolimits}                                    
\def\FijO{\FF^{\,i\to\bullet j}}                                
\def\FOij{\FF^{\,i\bullet\to j}}                                
\def\FtoO{\FF^{\,\to\bullet}}                                   
\def\FOto{\FF^{\bullet\to}}                                     
\def\cdc{,\ldots,}                                              
\def\1n{1\cdc n}                                                
\def\eq#1{\begin{equation}#1\end{equation}}                     
\def\eqs*#1{\begin{eqnarray*}#1\end{eqnarray*}}                 
\def\eqss#1{\begin{eqnarray}#1\end{eqnarray}}                   
\newtheorem{thm}{Theorem}{\bfseries}{\itshape}                  
\newtheorem{prop}{Proposition}{\bfseries}{\itshape}             
\newtheorem{corol}{Corollary}{\bfseries}{\itshape}              
\newtheorem{defin}{Definition}{\bfseries}{\upshape}             
\newtheorem{lemma}{Lemma}{\bfseries}{\itshape}                  
\def\proof{{\noindent\bf Proof. }}                              
\def\R{{\mathbb R}}                                             
\def\ovec{\overrightarrow}                                      
\def\T{{\rm\scriptscriptstyle T}}
\newtheorem{remark}{Remark}{\bfseries}{\upshape}                

\providecommand{\url}[1]{#1}
\csname url@samestyle\endcsname

\author{Pavel Chebotarev\\
       {\normalsize Institute of Control Sciences of the Russian Academy of Sciences}\\
       {\normalsize 65 Profsoyuznaya Street, Moscow 117997, Russia}\\
       {\small\tt chv@member.ams.org}
       }

\title{\vspace{-2.2em}The graph bottleneck identity}
\date{}

\begin{document}
\maketitle

\vspace{-2em}
\begin{abstract}
A matrix $S=(s_{ij})\in\R^{n\times n}$ is said to determine a \emph{transitional measure\/} for a digraph $\G$ on $n$ vertices if for all
$i,j,k\in\{\1n\},$ the \emph{transition inequality\/} $s_{ij}\,s_{jk}\le s_{ik}\,s_{jj}$ holds and reduces to the equality (called the \emph{graph bottleneck identity}) if and only if every path in $\G$ from $i$ to $k$ contains~$j$. We show that every positive transitional measure produces a distance by means of a logarithmic transformation. Moreover, the resulting distance $d(\cdot,\cdot)$ is \emph{graph-geodetic}, that is, $d(i,j)+d(j,k)=d(i,k)$ holds if and only if every path in $\G$ connecting $i$ and $k$ contains~$j$. Five types of matrices that determine transitional measures for a digraph are considered, namely, the matrices of path weights, connection reliabilities, route weights, and the weights of in-forests and out-forests.


\end{abstract}

\section{Introduction}

Two interesting properties of several well-known proximity/similarity measures $s(i,j)=s_{ij}$ for digraph vertices are that $s_{ij}\,s_{jk}\le s_{ik}\,s_{jj}$ and that $s_{ij}\,s_{jk}=s_{ik}\,s_{jj}$ if and only if every path from $i$ to $k$ contains~$j$. We call these the \emph{transition inequality\/} and the \emph{graph bottleneck identity}, respectively. For the path accessibility with a sufficiently small parameter and also for the connection reliability, the route accessibility, and two versions of the directed forest accessibility, the foregoing properties are proved in Sections~\ref{s_sii1} and~\ref{s_main} below.
In Sections~\ref{s_build} and~\ref{s_bottle}, we show that every positive-valued function with the above properties (we call such functions \emph{transitional measures}) gives rise to a graph-geodetic (i.e., such that $d(i,j)+d(j,k)=d(i,k)$ if and only if every path connecting $i$ and $k$ contains~$j$) logarithmic metric. As a synonym of metric, we use the term \emph{distance}, i.e., a distance is assumed to satisfy the triangle inequality. Graph-geodetic distances, in particular, are useful because they enable one to instantly check whether there are paths connecting $i$ and $k$ and not passing through~$j$ for any vertices $i$, $j$, and~$k$. Moreover, they have interesting mathematical properties.
In the rest of this section, we introduce some graph-theoretic notation and basic results mainly used in Sections~\ref{s_sii1} and~\ref{s_main}.

Let $\G$ be a weighted directed multigraph (in what follows, for brevity, a ``\emph{digraph\/} $\G$'') with vertex set $V=V(\G)=\{\1n\}$, $n>1$. Assume that $\G$ has no loops.
For $i,j\in V$, let $n_{ij}\in\{0,1,\ldots\}$ be the number of arcs emanating from $i$ to $j$ in~$\G$; for every $p\in\{\1n_{ij}\}$, let $w_{ij}^p>0$ be the weight of the $p$th arc directed from $i$ to $j$ in~$\G$; let $w_{ij}=\sum_{p=1}^{n_{ij}}w_{ij}^p$ (if $n_{ij}=0$, we set $w_{ij}=0$) and $W=(w_{ij})_{n\times n}$. $W$ is the \emph{matrix of total arc weights}. The \emph{outdegree\/} and \emph{indegree\/} of vertex $i$ are $\od(i)=\sum_{j=1}^nn_{ij}$ and $\id(i)=\sum_{j=1}^nn_{ji}$, respectively.

\smallskip
By the weight of a digraph $H$, $w(H)$, we mean the product of the weights of all its arcs. If $H$ has no arcs, then $w(H)=1$. The weight of a finite or denumerable set ${\cal S}$, $w({\cal S})$, is the sum of the weights of the elements in~${\cal S}$; the weight of the empty set is zero. If ${\cal S}$ is finite and contains digraphs whose arc weights are unity (i.e., the digraphs in ${\cal S}$ are actually unweighted), then $w({\cal S})$ is equal to the cardinality of~${\cal S}$.

\smallskip
For $v_0,v_k\in V(\G),$ a $v_0\to v_k$ \emph{path\/} in $\G$ is an alternating sequence of vertices and arcs $v_0,a_1,v_1\cdc a_k,v_k$ where all vertices are distinct and each $a_i$ is a $v_{i-1}\!\to\!v_i$ arc. The unique $v_0\to v_0$ path is the ``sequence''\,$v_0$ having no arcs.
The \emph{length\/} of a path is the number $k$ of its arcs. The \emph{weight\/} of a path is the product of the weights of its arcs.
The weight of a $v_0\to v_0$ path is~1.
A~digraph is \emph{strong\/} (or \emph{strongly connected}) if for every vertices $v$ and $v',$ it has a $v\to v'$ path. A~digraph is \emph{weakly connected\/} if the corresponding undirected graph is connected.

\smallskip
A \emph{converging tree\/} is a weakly connected weighted digraph in which one vertex, called the {\it root}, has outdegree zero and the remaining vertices have outdegree one.
A~{\it converging forest\/} is a weighted digraph all of whose weakly connected components are converging trees. The roots of these trees are referred to as the roots of the converging forest. A spanning converging forest of $\G$ is called an {\em in-forest\/} of~$\G$.

\smallskip
For a fixed digraph $\G$, by $\FtoO$ and $\FijO$ we denote the set of all in-forests of $\G$ and the set of all in-forests of $\G$ that have vertex $i$ belonging to a tree rooted at~$j$, respectively. Let
$$
f=w(\FtoO)
$$
and
\eq{
\label{e_fij} %
f_{ij}=w(\FijO),\quad i,j\in V;
}
$F=(f_{ij})_{n\times n}$ is called the \emph{matrix of in-forests of~\,$\G$}.

\smallskip
Let $L=(\ell_{ij})$ be the Laplacian matrix of $\G$, i.e., for $\;i,j=\1n$,
\eq{
\label{e_Lapl}
\ell_{ij}=
        \begin{cases}
        -w_{ij},               &j\ne i,\\
         \suml_{k\ne i}w_{ik}, &j  = i.
        \end{cases}
}

Consider the matrix
\eq{
\label{e_Q}
Q=(q_{ij})=(I+L)^{-1},
}
where $I$ is the identity matrix. By the matrix forest theorem \cite{CheSha97,CheAga02ap} (``undirected'' versions of this theorem can be found in~\cite{CheSha95,Merris97}), for any digraph $\G$, $Q$ does exist and
\eq{
\label{e_mft}
q_{ij}=\frac{f_{ij}}{f},\quad i,j=\1n.
}
Therefore, $F=fQ=f\cd(I+L)^{-1}$. The matrix $Q$ can be considered as a proximity (similarity) matrix of~$\G$ \cite{CheSha97,Che08DAM}; it has a random walk interpretation \cite[Section~4]{Che08DAM}; in the case of undirected graphs, it is also called the \emph{regularized Laplacian kernel} (cf.~\cite{SmolaKondor03}).

\smallskip
In Sections~\ref{s_sii1} and~\ref{s_main}, we show that the values $f_{ij}$ and several other proximity indices satisfy the transition inequality and the graph bottleneck identity. Some general implications of these properties (mainly relating to the construction of graph distances) are studied in Sections~\ref{s_trans}, \ref{s_build}, and~\ref{s_bottle}. The results obtained have undirected counterparts; one of them is presented in Section~\ref{s_undir}. In \cite{Che08dist}, the approach of this paper is used to fill the gap between the shortest path distance and the resistance distance for undirected graphs.

\section{Transitional measures and the graph bottleneck\\ identity}
\label{s_trans}

We say that a matrix $S\!=\!(s_{ij})\!\in\!\R^{n\times n}$ satisfies the \emph{transition inequality\/} if for all $1\le i,j,k\le n,$
\eq{
\label{e_trans}
s_{ij}\,s_{jk}\le s_{ik}\,s_{jj}.
}

\begin{lemma}
\label{l_trans}
If $S\!=\!(s_{ij})\!\in\!\R^{n\times n}$ satisfies the transition inequality$,$ then
for all\/ $1\le i,j\le n,$
\eq{
\label{e_KB}
s_{ij}\,s_{ji}\le s_{ii}\,s_{jj}.
}
\end{lemma}

\proof
This is immediate by setting $k=i$ \/in~\eqref{e_trans}.
\qed

\begin{remark}
{\rm
Inequality~\eqref{e_KB} bears a close analogy to the Cauchy--Bunyakovsky--Schwarz inequality. Therefore, if $S$ is symmetric, has positive diagonal, and satisfies \eqref{e_trans}, then it can be treated as a matrix of variances and covariances or a Gram matrix. As a result, say, $\arccos\tfrac{s_{ij}}{\sqrt{s_{ii}\,s_{jj}}}$ can be considered as the angle between the objects represented by~$i$ and~$j$, which is suitable for scaling purposes; see also \cite[Section~7.9]{BorgGroenen05}. 
At last the transition inequality is a multiplicative analogue of the \emph{triangle inequality for proximities\/} \cite{CheSha97,CheSha98} also called the ``unrooted correlation triangle inequality''~\cite{DezaLaurent97}.
}
\end{remark}

Furthermore, we say that a matrix $S\!=\!(s_{ij})\!\in\!\R^{n\times n}$ satisfies the \emph{graph bottleneck identity\-\/} w.r.t.\ a digraph $\G$ (an undirected multigraph $G$) with vertex set $V=\{\1n\}$ if for all ${1\le i,j,k\le n,}$
\eq{
\label{e_boe}
s_{ij}\,s_{jk}=s_{ik}\,s_{jj}
}
holds if and only if all directed paths in $\G$ (all paths in $G$) from $i$ to $k$ contain~$j$.
\medskip

Eq.~\eqref{e_boe} is referred to as the graph bottleneck identity because it pertains to the case where $j$ is a kind of a bottleneck (or a cut point) for the $i\to k$ paths: the removal of $j$ disconnects $k$ from~$i.$

To shorten the terminology, we give the following definition.

\begin{defin}
\label{def_trme}
{\rm
Given a digraph $\G$ with vertex set $V=\{\1n\},$ suppose that a matrix $S=(s_{ij})_{n\times n}$ satisfies the transition inequality~\eqref{e_trans} and the graph bottleneck identity~\eqref{e_boe} w.r.t.~$\G$.\ Then we say that $S$ \emph{determines the transitional measure $s(i,j)\!=\!s_{ij},$ $i,j\!\in\! V,$ for $\G$}.
}
\end{defin}

For undirected graphs, the notion of transitional measure is defined similarly.
It will be shown in Sections \ref{s_sii1} and \ref{s_main} that several popular graph proximity measures are transitional.

\begin{lemma}
\label{l_bottl}
If $S\!=\!(s_{ij})\!\in\!\R^{n\times n}$ determines a transitional measure for some digraph $\G,$ then\footnote{Inequality~\eqref{e_KBs} also holds for every matrix $S$ that, with no relation to graphs, obeys the \emph{strengthened transition inequality}, which is \eqref{e_trans} turning into the strict form whenever $k=i$ and $j\ne i$. It follows from the proof of Theorem~\ref{th_dist} that if such a matrix has positive off-diagonal entries, then it 
produces a distance by means of \eqref{Ha} and~\eqref{Da}.}
for all\/ $1\le i,j\le n$ such that $j\ne i,$
\eq{
\label{e_KBs}
s_{ij}\,s_{ji}<s_{ii}\,s_{jj}.
}
\end{lemma}

\proof
Setting $k=i$ in~\eqref{e_trans} and taking into account that there is a path of length 0 from $i$ to $k=i$ that does not contain $j\ne i$\/ we conclude that the transition inequality and the graph bottleneck identity yield~\eqref{e_KBs}.
\qed

\medskip
The main object of our interest in this paper is the distances constructed on the basis of transitional measures.

\section{Logarithmic distances built on the basis of\\ the transition inequality}
\label{s_build}

If a matrix $S$ satisfies the transition inequality \eqref{e_trans} and its off-diagonal entries are positive, then all the entries of $S$ are positive.
In this case, define the matrix
\eq{
\label{Ha}
H=\ovec{\ln S},
}
where
$\ovec{\varphi(S)}$ stands for elementwise operations, i.e., operations applied to each entry of $S$ separately. Consider the matrix
\eq{
\label{Da}
D=\tfrac{1}{2}(h{\bm1}^\T+\bm1 h^\T-H-H^\T),
}
where $h$ is the column vector containing the diagonal entries of $H$, ${\bm1}$ is the column of $n$ ones, and
$H^\T$, $h^\T$, and ${\bm1}^\T$ are the transposes of $H$, $h$, and ${\bm1}$.
An alternative form of \eqref{Da} is $D=(U+U^\T)/2$, where $U=h{\bm1}^\T-H$, and the elementwise form is $d_{ij}=\frac{1}{2}(h_{ii}+h_{jj}-h_{ij}-h_{ji}),\,$ $i,j=\1n,$ where $H=(h_{ij})$ and $D=(d_{ij}).$ This is a standard transformation used to obtain a distance from a proximity measure (cf.\ the inverse covariance mapping in~\cite{DezaLaurent97} and \cite[Section~12.1]{BorgGroenen05}).

\begin{thm}
\label{th_dist}
If $S\!=\!(s_{ij})_{n\times n}$ determines a transitional measure for some digraph\/ $\G$ and has positive off-diagonal entries$,$ then
$D=(d_{ij})_{n\times n}$ defined by \eqref{Ha} and \eqref{Da} is a matrix of distances on\/~$\{\1n\}.$
\end{thm}

Before proving Theorem~\ref{th_dist} we give an expression for the entries of $D$. Eqs.\ \eqref{Ha} and \eqref{Da} for every $i,j=\1n$ imply
\eqss{
\label{p_d1}
d_{ij}
=\tfrac{1}{2}(h_{ii}+h_{jj}-h_{ij}-h_{ji})
=\tfrac{1}{2}(\ln s_{ii}+\ln s_{jj}-\ln s_{ij}-\ln s_{ji})
=\tfrac{1}{2}\ln\!\frac{s_{ii}\,s_{jj}}{s_{ij}\,s_{ji}}.
}

{\noindent\bf Proof of Theorem~\ref{th_dist}.}
The proof amounts to showing that for all $i,j,k\in \{\1n\}$:

(i) $d_{ij}=0$ if and only if $i=j\,$ and

(ii) $d_{ij}+d_{jk}-d_{ki}\ge0$ (triangle inequality).

\smallskip
Indeed, the symmetry and non-negativity of $D$, which are sometimes considered as part of the definition of distance, follow from (i) and~(ii).
Since $S$ \/has positive off-diagonal entries, the transition inequality implies the positivity of~$S$.

To prove (i), note that if $i=j,$ then by~\eqref{Da}, $d_{ij}=0$. Conversely, if $d_{ij}=0$, then by~\eqref{p_d1}, $s_{ii}\,s_{jj}=s_{ij}\,s_{ji}$ holds, which, by Lemma~\ref{l_bottl}, implies that $i=j$.

\smallskip
To prove (ii), observe that by \eqref{Ha}, \eqref{Da}, and the transition inequality \eqref{e_trans},
\eqss{
\label{e_trine}
d_{ij}+d_{jk}-d_{ki}
&=&\tfrac{1}{2}(h_{ii}+h_{jj}+h_{jj}+h_{kk}-h_{kk}-h_{ii}\nonumber\\
&&\,-h_{ij}-h_{ji}-h_{jk}-h_{kj}+h_{ki}+h_{ik})\nonumber\\
&=&\tfrac{1}{2}\ln\Bigl(\frac{s_{jj}\,s_{ik}}{s_{ij}\,s_{jk}}\!\cdot\!\frac{s_{jj}\,s_{ki}}{s_{kj}\,s_{ji}}\Bigr)\ge0
}
holds. This completes the proof.
\qed

\medskip
Based on Theorem~\ref{th_dist}, we give the following definition.

\begin{defin}
\label{def_dist}
{\rm
Suppose that $S\!=\!(s_{ij})_{n\times n}$ has positive off-diagonal entries and determines a transitional measure for some digraph~$\G$.
The \emph{logarithmic distance corresponding to $S$\/} is the function $d\!:\{\1n\}^2\to\R$ such that $d(i,j)=d_{ij},\,$ $i,j=\1n,$ where $D=(d_{ij})$ is defined by \eqref{Ha} and~\eqref{Da}.
}
\end{defin}

In Section~\ref{s_bottle}, it is shown that every distance of  this kind is graph-geodetic.

\section{The graph bottleneck identity implies the geodetic\\ property of the logarithmic distance}
\label{s_bottle}

\begin{defin}
\label{d_g-d}
{\rm
For a multidigraph $\G$ (a multigraph $G$) with vertex set $V,$ a function $d\!:V\!\times\!V\to\R$ is called \emph{graph-geodetic\/} provided that
$d(i,j)+d(j,k)=d(i,k)$ holds if and only if every directed path in $\G$ connecting $i$ and $k$ in either direction (every path in $G$ connecting $i$ and~$k$) contains~$j$.}
\end{defin}

If $d(\cdot,\cdot)$ is a distance on digraph vertices, then the property of being graph-geodetic (this term is taken from~\cite{KleinZhu98}) is a natural condition of strengthening the triangle inequality to equality. Knowing a graph-geodetic distance enables one to instantly check whether $j$ ``separates'' $i$ and $k$ or not for any $i,j,k\in V(\G).$ The classical shortest path distance clearly possesses the ``if'' (but not the ``only if'') part of the graph-geodetic property; the ``if'' part of this property for the resistance distance was proved in~\cite{KleinRandic93}. The ordinary distance in a Euclidean space satisfies a similar condition resulting from substituting ``line segment'' for ``path in $G$.''

\begin{thm}
\label{th_botadd}
Suppose that $S\!=\!(s_{ij})_{n\times n}$ has positive off-diagonal entries and determines a transitional measure for some digraph\/ $\G$.
Then the logarithmic distance corresponding to $S$ is graph-geodetic for\/~$\G$.
\end{thm}

\proof
Using \eqref{e_trine} and the transition inequality we conclude that
$d_{ij}+d_{jk}=d_{ki}$ is true if and only if $\tfrac{s_{jj}\,s_{ik}}{s_{ij}\,s_{jk}}=\tfrac{s_{jj}\,s_{ki}}{s_{kj}\,s_{ji}}=1$. In turn, by the graph bottleneck identity, this holds if and only if every path in $\G$ connecting $i$ and $k$ in either direction contains~$j.$ Thus, by Definition~\ref{d_g-d}, the logarithmic distance $d(i,j)=d_{ij}$ ($i,j=\1n$) corresponding to $S$ is graph-geodetic for~$\G$.
\qed

\medskip
Graph-geodetic functions have many interesting properties. One of them, as mentioned in \cite{KleinRandic93}, is a simple connection (such as that obtained in \cite{GrahamHoffmanHosoya77}) between the cofactors and the determinant of $\G$'s distance matrix and those of the maximal blocks of $\G$ that have no cut points. Another property is the recursive Theorem~8 in~\cite{KleinZhu98}. The graph-geodetic distances are not Euclidean; however, by Blumenthal's ``Square-Root'' theorem, the corresponding ``square-rooted'' distances satisfy the 3-Euclidean condition (see, e.g.,~\cite{KleinZhu98}).

\smallskip
Obviously, it is \eqref{Ha} that guarantees the graph-geodetic property of the matrix $D$ obtained by means of~\eqref{Da} from a transitional measure. 
If $H=S,$ then this property is not secured and a sufficient condition of $D$'s being a distance matrix is provided by the following proposition. 

\begin{prop}
\label{p_noLog}
Suppose that $S=(s_{ij})_{n\times n}$ satisfies the transition inequality~\eqref{e_trans} and
\eq{
\label{e_dido}
s_{jj}>\min(s_{ij},s_{ji}),\;\;s_{jj}\ge\max(s_{ij},s_{ji}),\;\;\text{and}\;\;s_{jj}>0\;\;\text{for all}\;\;\, i,j=\1n,\; j\ne i.
}
Then $D$ defined by \eqref{Da} with $H=S$ is a matrix of distances.
\end{prop}

\proof
Assuming that~\eqref{e_trans} and \eqref{e_dido} are satisfied we prove that
(i)~$d_{ij}=0$ if and only if $i=j\,$ and
(ii)~$d_{ij}+d_{jk}-d_{ki}\ge0$
for all $i,j,k=\1n.$
Since by \eqref{Da},
\eqss{
\nonumber
d_{ij}
&=&\tfrac{1}{2}(s_{ii}+s_{jj}-s_{ij}-s_{ji})\;\,\text{and}\\
\label{e_ddnoL}
d_{ij}+d_{jk}-d_{ki}
&=&\tfrac{1}{2}((s_{jj}+s_{ik}-s_{ij}-s_{jk})+(s_{jj}+s_{ki}-s_{kj}-s_{ji}))
}
hold, $(j=i)\Rightarrow (d_{ij}=0)$ is immediate and $(j\ne i)\Rightarrow (d_{ij}\ne0)$ follows from \eqref{e_dido}. Furthermore, since by \eqref{e_dido}, $s_{jj}>0,$ \eqref{e_trans} implies that $s_{ik}\ge s_{ij}s_{jk}s_{jj}^{-1}$ and $s_{ki}\ge s_{kj}s_{ji}s_{jj}^{-1}$, therefore, by \eqref{e_ddnoL} and~\eqref{e_dido},
\eqs*{\qquad\,
d_{ij}+d_{jk}-d_{ki}
&\ge&\tfrac{1}{2} ((s_{jj}+s_{ij}\,s_{jk}\,s_{jj}^{-1}-s_{ij}-s_{jk})+
                   (s_{jj}+s_{kj}\,s_{ji}\,s_{jj}^{-1}-s_{kj}-s_{ji}))\\
& = &\tfrac{1}{2}(((s_{jj}-s_{ij})(s_{jj}-s_{jk})
                 + (s_{jj}-s_{ji})(s_{jj}-s_{kj}))s_{jj}^{-1})\ge0.\qquad\;\qed
}

\smallskip
In Sections~\ref{s_sii1} and~\ref{s_main}, we show that several well-known graph proximity measures are transitional.

\section{Two transitional measures with unit diagonal}
\label{s_sii1}

In this section, we consider two instances of transitional measures. With relation to the graph bottleneck identity, they represent a very special case 
in which for every $i\in V,\,$ $s_{ii}=1.$

\subsection{The path $\tau$-accessibility}
The \emph{path $\tau$-accessibility\/} of $j$ from $i$ in $\G$ is the total \emph{$\tau$-weight\/} of all paths from $i$ to~$j$:
\eq{
\label{e_path}
s_{ij}=w_\tau(\PP^{ij})=\sum_{P_{ij}\in\PP^{ij}}w_\tau(P_{ij}),
}
where $\PP^{ij}$ is the set of all $i\to j$ paths in~$\G,$
\eqs*{
w_\tau(P_{ij})=\tau^{l(P_{ij})}w(P_{ij}),
}
$l(P_{ij})$ and $w(P_{ij})$ are the length and the weight of $P_{ij},$ and $\tau>0.$

By definition, for every $i\in V,$ the unique ``path from $i$ to $i$'' is the path of length $0$ whose weight is unity, whence $s_{ii}=1,\;i=\1n$.

\begin{thm}
\label{th_paths}
For any digraph $\G,$ there exists $\tau_0>0$ such that for every\/ $\tau\in(0,\tau_0),$ $S=(s_{ij})$ defined by~\eqref{e_path}
determines a transitional measure for~$\G.$
\end{thm}

\proof
For arbitrary $i,j,k\in V,$ $P_{ij}\in\PP^{ij},$ and $P_{jk}\in\PP^{jk}$, let $v$ be the first (along $P_{ij}$) vertex of $P_{ij}$ that belongs to $P_{jk}$. Then combining the $i\to v$ subpath of $P_{ij}$ with the $v\to k$ subpath of $P_{jk}$ we obtain a well-defined path $P_{ik}\in\PP^{ik}$ whose $\tau$-weight is no less than $w_\tau(P_{ij})\!\cdot\!w_\tau(P_{jk})$ for each sufficiently small $\tau>0$. If this $P_{ik}$ contains $j$ (i.e., $v=j$), then
\eq{
\label{e_Pik1} 
w_\tau(P_{ik})=w_\tau(P_{ij})\,w_\tau(P_{jk})
}
for every $\tau>0.$ Otherwise, if a fixed $P_{ik}$ does not contain $j,$ then a $\tau_0(P_{ik},j)>0$ can be chosen in such a way that
\eq{
\label{e_Pik2} 
w_\tau(P_{ik})>\sum_{(P_{ij},\,P_{jk})\to P_{ik}}\!\!\!w_\tau(P_{ij})\,w_\tau(P_{jk})
}
for all $0<\tau<\tau_0(P_{ik},j),$ where the sum is taken over all $P_{ij}\in\PP^{ij}$ and $P_{jk}\in\PP^{jk}$ such that combining the $i\to v$ subpath of $P_{ij}$ with the $v\to k$ subpath of $P_{jk}$ produces the fixed $P_{ik}$ (which is denoted by $(P_{ij},\,P_{jk})\to P_{ik}$). Let $\tau_0=\min_{i,\,j,\,k\in V,\,P_{ik}\in\PP^{i\bar\jmath k}} \{\tau_0(P_{ik},j)\},$ where $\PP^{i\bar\jmath k}$ is the set of all $i\to k$ paths in $\G$ that do not contain~$j.$ Thus, if $0<\tau<\tau_0,$ then \eqref{e_Pik2} holds for all $P_{ik}\in\PP^{i\bar\jmath k}$ and \eqref{e_Pik1} holds for all $P_{ik}\in\PP^{ik}\smallsetminus\PP^{i\bar\jmath k}$. Consequently, for any $\tau\in(0,\tau_0)$ and any $i,j,k\in V,$
\eqs*{
s_{ik}\,s_{jj}
&=&s_{ik}
=  \sum_{ P_{ik}\in\PP^{ik}}             w_\tau(P_{ik})
\ge\sum_{ P_{ik}\in\PP^{ik}}
   \sum_{(P_{ij},\,P_{jk})\to P_{ik}}\!\!w_\tau(P_{ij})\,w_\tau(P_{jk})\\
&=&\sum_{ P_{ij}\in\PP^{ij}}         \!\!w_\tau(P_{ij})
   \sum_{ P_{jk}\in\PP^{jk}}         \!\!w_\tau(P_{jk})
=  s_{ij}\,s_{jk},
}
with the equality if and only if every $i\to k$ path contains~$j$. The transition inequality and the graph bottleneck identity follow.
\qed

\subsection{Connection reliability}
\label{s_conRel}
Consider a digraph $\G$ with arc weights $w_{ij}^p\in(0,1]$ interpreted as the intactness probabilities of the arcs. Define $p_{ij}$ to be the $i\to j$ \emph{connection reliability\/}, i.e., the probability that at least one path from $i$ to $j$ remains intact, 
provided that the arc failures are independent. Let $P=(p_{ij})$ be the matrix of connection reliabilities for all pairs of vertices. For every $j\in V,\,$ $p_{jj}=1$, because the $j\to j$ path of length~$0$ is always intact.

The connection reliabilities can be represented as follows (see, e.g., \cite[p.~10]{Shier}):

\vspace{-.9em}
\eq{
\label{e_Pways}
p_{ij}    =\suml_k        \!\Pr(P_k)
          -\suml_{k<t  }  \!\Pr(P_k P_t)
          +\suml_{k<t<l}\!\!\Pr(P_k P_t P_l)-\ldots
          +(-1)^{m-1}       \Pr(P_1 P_2\cdots P_m),
}

\vspace{-.3em}\noindent
where $P_1, P_2\cdc P_m$ are all $i\to j$ paths in $\G$, $\Pr(P_k)=w(P_k),$ $\Pr(P_kP_t)=w(P_k\cup P_t)$, $P_k\cup P_t$ is the subdigraph of $\G$ containing those arcs that belong to $P_k$ or $P_t$, and so forth. By virtue of \eqref{e_Pways}, connection reliability is a modification of path accessibility that takes into account the degree of overlap for various paths between vertices.

\vspace{-.2em}
\begin{thm}
\label{th_reli}
For any digraph $\G$ with arc weights $w_{ij}^p\in(0,1],$ the matrix $P=(p_{ij})$ of connection reliabilities determines a transitional measure for~$\G.$
\end{thm}

\vspace{-.2em}
\proof
Let $E_{ij}$ be the event that at least one path connecting $i$ to $j$ remains intact. Then, since $E_{ij}\wedge E_{jk}\Rightarrow\!E_{ik}$, by the independence assumption we have

\vspace{-1.4em}
\eqs*{
p_{ik}\,p_{jj}=p_{ik}=\Pr(E_{ik})\ge\Pr(E_{ij})\Pr(E_{jk})=p_{ij}\,p_{jk}
}

\vspace{-.3em}\noindent
with the equality if and only if every path from $i$ to $k$ contains~$j$.
\qed

\vspace{-.2em}
\begin{corol}[{of Theorems~\ref{th_botadd}, \ref{th_paths}, and~\ref{th_reli}}]
\label{co_paths}
For any strong digraph $\G,$ the logarithmic distances corresponding to the matrix $S=(s_{ij})$ defined by~\eqref{e_path} with a sufficiently small $\tau$ and to the matrix $P=(p_{ij})$ of connection reliabilities $($whenever $w_{ij}^p\in(0,1])$ are graph-geodetic for~$\G$.
\end{corol}

\vspace{-.2em}
\proof
Since for a strong digraph $\G,$ the matrices $S$ and $P$ have positive off-diagonal entries, the desired statements follow from Theorems~\ref{th_paths}, \ref{th_reli}, and~\ref{th_botadd}.
\qed

\smallskip
The next section is devoted to the transitional measures in which the diagonal elements $s(i,i)$ measure the (relative) strength of connections of every vertex to itself.

\vspace{-.2em}
\section{The matrices of spanning forests and routes provide\\ transitional measures}
\label{s_main}

The following theorem is the main technical result of this paper.

\vspace{-.3em}
\begin{thm}
\label{th_fores}
For any digraph $\G,$ the matrix of in-forests $F=(f_{ij})$ defined by~\eqref{e_fij} determines a transitional measure for~$\G.$
\end{thm}

\vspace{-.3em}
There seems to be no easy way to construct a direct bijective proof of Theorem~\ref{th_fores} (such as the proofs of Theorems~\ref{th_paths} and~\ref{t_Routes}). So we present an indirect proof relying on Proposition~\ref{p_RoutFo} and Theorem~\ref{t_Routes} given below. We will use the following construction.

For a fixed digraph $\G$, let us choose an arbitrary $\e>0$ such that
\eq{
\label{e_rest}
\e\maxl_{1\le i\le n}\ell_{ii}<1,
}
where $L=(\ell_{ij})$ is the Laplacian matrix of $\G$, whose diagonal entries are always non-negative (see~\eqref{e_Lapl}).
It is easy to see that the matrix
\eq{
\label{e_Pm}
P=(p_{ij})=I-\e L
}
(not to be confused with the matrix $P$ of Section~\ref{s_conRel}) is row stochastic: $0\le p_{ij}\le1$ and $\sum_{k=1}^n p_{ik}=1$, $\;i,j=\1n$.
\smallskip\smallskip

Denote by $\G^{\circlearrowright}$ a weighted multidigraph with loops whose matrix of total arc weights is
\eq{
\label{e_glu}
W(\G^{\circlearrowright})=(1+\e)^{-1}P.
}
$\G^{\circlearrowright}$ can be constructed as follows: every vertex $i$ of $\G^{\circlearrowright}$ gets a loop with weight $(1+\e)^{-1}p_{ii}$; the remaining arcs of $\G^{\circlearrowright}$ are the same as in $\G$, their weights being equal to the corresponding weights in~$\G$ multiplied by $(1+\e)^{-1}\e$.

Recall that a $v_0\to v_k$ \emph{route\/} (also called a \emph{walk}) in a multidigraph with loops is an arbitrary alternating sequence of vertices and arcs $v_0,a_1,v_1\cdc a_k,v_k$ where each $a_i$ is a $v_{i-1}\!\to\!v_i$ arc. 
The \emph{length\/} of a route is the number $k$ of its arcs (including loops). The \emph{weight\/} of a route is the product of the $k$ weights of its arcs (including repeated arcs). By definition, for every vertex $v_0$, there is a $v_0\to v_0$ route $v_0$ with length $0$ and weight~1.
\smallskip

Let $r_{ij}$ be the weight of the set $\RR^{ij}$ of all $i\to j$ routes in $\G^{\circlearrowright}$, provided that this weight is finite (note that in the presence of loops $\RR^{ij}$ is infinite whenever $j$ is reachable from~$i$). $R=(r_{ij})_{n\times n}$ will denote the \emph{matrix of the route weights}.

\begin{prop}
\label{p_RoutFo}
For any digraph $\G$ and any $\e\!>\!0$ that satisfies \eqref{e_rest}$,$ the matrix $R$ of the route weights in $\G^{\circlearrowright}\!$ is finite and positively proportional to the matrix $F$\! of in-forests of\/~$\G$.
\end{prop}

\proof
By \eqref{e_glu}, for each $k=0,1,2,\ldots,$ the matrix of the weights of $k$-length routes in $\G^{\circlearrowright}$ is $((1+\e)^{-1}P)^k$. Therefore, the matrix $R$, whenever it exists, can be represented as follows:
\eq{
\label{e_RfrP}
R=\suml_{k=0}^\infty((1+\e)^{-1}P)^k.
}
Since the spectral radius of $P$ is 1 and $0<(1+\e)^{-1}<1$, the series in \eqref{e_RfrP} converges to a finite matrix\footnote{On counting routes, see also~\cite{Kasteleyn67,ChelnokovZefirova09}. Related finite topological representations that involve paths are obtained in~\cite{Ponstein66}. For some connections with matroid theory, we refer to~\cite{Schrijver78}.}, 
therefore \eqref{e_RfrP}, \eqref{e_Pm}, \eqref{e_Q}, and \eqref{e_mft} imply
\eqs*{
R
&=&(I-(1+\e)^{-1}P)^{-1}
 = \left(I-(1+\e)^{-1}\!\left(I-\e L\right)\!\right)^{-1}\\
&=&\left(\frac{\e}{1+\e}     (I+   L)      \!\right)^{-1}
 = \left(1+\e^{-1}\right)Q
 = \left(1+\e^{-1}\right)f^{-1}F,
}
which completes the proof.\qed

\begin{thm}
\label{t_Routes}
For any weighted multidigraph allowing loops$,$
if the matrix ${R=(r_{ij})_{n\times n}}$ of route weights is finite$,$ then $R$ determines a transitional measure for this multidigraph.
\end{thm}

\proof
Suppose that $R$ is finite.
Let $\RR^{ij(1)}$ be the set of all $i\to j$ routes that contain only one occurrence of~$j$. Let $r_{ij(1)}=w(\RR^{ij(1)})$. Then every $i\to j$ route $\rr_{ij}\in \RR^{ij}$ can be uniquely decomposed into a route $\rr_{ij(1)}\in\RR^{ij(1)}$ and a route $\rr_{jj}\in\RR^{jj}$ (if $\rr_{ij}\in \RR^{ij(1)},$ then $\rr_{ij}$ is decomposed into itself and the $j\to j$ route of length~0). And vice versa, linking an arbitrary route $\rr_{ij(1)}\in\RR^{ij(1)}$ with an arbitrary $\rr_{jj}\in\RR^{jj}$ results in a well-defined route $\rr_{ij}\in\RR^{ij}$. This induces a natural bijection between $\RR^{ij}$ and $\RR^{ij(1)}\times\RR^{jj}$. Therefore
\eq{
\label{e_rij}
r_{ij}=r_{ij(1)}\,r\!_{jj}.
}

Let $\RR^{ijk}$ and $\RR^{i\bar\jmath k}$ be the sets of all $i\to k$ routes that contain and do not contain $j$, respectively. Then $\RR^{ik}=\RR^{ijk}\cup\RR^{i\bar\jmath k}$ and $\RR^{ijk}\cap\RR^{i\bar\jmath k}=\varnothing$, consequently,
\eq{
\label{e_rik}
r_{ik}=r_{ijk}+r_{i\bar\jmath k},
}
where $r_{ijk}=w(\RR^{ijk})$ and $r_{i\bar\jmath k}=w(\RR^{i\bar\jmath k})$.
\smallskip

Furthermore, by the argument similar to that justifying \eqref{e_rij} one has
\eq{
\label{e_rijk}
r_{ijk}=r_{ij(1)}\,r\!_{jk}.
}

Combining \eqref{e_rik}, \eqref{e_rijk}, and \eqref{e_rij} yields
$$
 r_{ik}\,r\!_{jj}
=(r_{ijk}+r_{i\bar\jmath k})\,r\!_{jj}
= r_{ij(1)}\,r\!_{jk}\,r\!_{jj}+r_{i\bar\jmath k}\,r\!_{jj}
= r_{ij}\,r\!_{jk}+r_{i\bar\jmath k}\,r\!_{jj}\geq r_{ij}\,r\!_{jk},
$$
with the equality if and only if $\RR^{i\bar\jmath k}=\varnothing$ (since $r_{jj}\ge1$).
The transition inequality and the graph bottleneck identity follow.
\qed
\bigskip

{\noindent\bf Proof of Theorem~\ref{th_fores}.}
Theorem~\ref{th_fores} is concluded by combining Proposition~\ref{p_RoutFo} and Theorem~\ref{t_Routes}.
\qed

\begin{corol}[{of Theorems~\ref{th_botadd}, \ref{th_fores}, and~\ref{t_Routes}}]
\label{co_fordi}
$1.$ For any strong digraph $\G,$ the logarithmic distance corresponding to the matrix of in-forests $F=(f_{ij})$ defined by~\eqref{e_fij} is graph-geodetic for~$\G$.

$2.$ For any strong weighted multidigraph allowing loops$,$ if the matrix 
${R=(r_{ij})_{n\times n}}$ is finite$,$ then the logarithmic distance corresponding to $R$ is graph-geodetic for this multidigraph.
\end{corol}

\proof
In view of Theorem~\ref{th_botadd}, the desired statements follow from Theorems~\ref{th_fores} and~\ref{t_Routes}.
\qed

\begin{remark}
{\rm
In Theorem~\ref{th_fores} and Corollary~\ref{co_fordi}, the matrix of in-forests $F=(f_{ij})$ can be replaced by the \emph{matrix $F'=(f'_{ij})$ of out-forests of\/}~$\G$. In greater detail, a spanning subdigraph $H$ of $\G$ is an \emph{out-forest\/} if every weak component of $H$ has one vertex of indegree zero (the root) and all other vertices of indegree one.
Consider the matrix
$
Q'=(q'_{ij})=(I+L')^{-1},
$
where $L'=(\ell'_{ij})$ is the \emph{column Laplacian\/} matrix~\cite[Section~2.2]{CheAga02ap} of~$\G$ whose entries are:
\eqs*{
\ell'_{ij}=
        \begin{cases}
        -w_{ij},               &j\ne i,\\
         \suml_{k\ne j}w_{kj}, &j  = i
        \end{cases}
}
(cf.\ \eqref{e_Lapl}--\eqref{e_Q}).
By the matrix forest theorem, $Q'$ does exist and
$
q'_{ij}={f'_{ij}}/{f'},\; i,j=\1n,
$
where $f'$ is the total weight of the out-forests in $\G$ ($f'=w(\FOto)$) and $f'_{ij}$ the total weight of out-forests having $j$ in a weak component rooted at~$i$ ($f'_{ij}=w(\FOij)$).

From these definitions it follows that $F'$ is the transposed matrix $F$ of the reverse digraph~$\G^{-1}$. Therefore, by Theorem~\ref{th_fores}, $F'$ determines a transitional measure for $\G$ and, in view of Theorem~\ref{th_botadd}, the corresponding logarithmic distance is graph-geodetic for~$\G$. It is worth noting that the logarithmic distances produced by $F$ and $F'$ are generally different.
}
\end{remark}

Finally, we touch upon the case of undirected graphs. This case is also considered in~\cite{Che08dist}.

\vspace{-.1em}
\section{On transitional measures for undirected graphs}
\label{s_undir}

For undirected multigraphs, the definitions of transitional measure and logarithmic distance are completely similar to Definitions~\ref{def_trme} and~\ref{def_dist}, and the above theorems have undirected counterparts. In this section, we present the least obvious result of this kind, which concerns spanning forests.

\begin{corol}[{of Theorem~\ref{th_fores}}]
\label{co_main}
Let $G$ be a connected weighted undirected multigraph and let $f_{ij},$ $i,j\in V(G),$ be the total weight of the spanning rooted forests of $G$ that have vertex $i$ belonging to a tree rooted at~$j$.
Then\/$:$

$1.$ The matrix $F=(f_{ij})$ determines a transitional measure for~$G;$

$2.$ The logarithmic distance corresponding to $F=(f_{ij})$ is graph-geodetic for~$G$.
\end{corol}

\proof
1. Consider the symmetric multidigraph $\G$ obtained from $G$ by replacing every edge by two opposite arcs carrying the weight of that edge. Then comparing the matrix forest theorems for directed and undirected graphs~\cite{CheSha97} yields $f_{ij}(G)=f_{ij}(\G),$ $i,j\in V(G)$. Observe that for every $i,j,k\in V(G),$ every path from $i$ to $k$ contains~$j\/$ if and only if so does every directed path from $i$ to $k$ in~$\G$. Therefore, by virtue of Theorem~\ref{th_fores}, $F=(f_{ij})$ determines a transitional measure for~$G.$
Item 2 follows from item 1 of Corollary~\ref{co_fordi}.
\qed

\vspace{-.1em}
\section*{Acknowledgements}

This work was partially supported by the RFBR Grant 09-07-00371 and the RAS Presidium Program ``Development of Network and Logical Control.''
The author is grateful to the referees for their comments.

\end{document}